\documentclass[11 pt]{article}

\usepackage{amssymb}

\def\nmark{\mbox{$\rm\bf\kern0.2em\rule{0.06em}{1.45ex}\kern-0.3em
N$}}
\def\dmark{\mbox{$\rm\bf\kern0.2em\rule{0.06em}{1.45ex}\kern-0.3em
D$}}
\def\cmark{\mbox{$\rm\bf\kern0.2em\rule{0.06em}{1.45ex}\kern-0.3em
C$}}
\def\rmark{\mbox{$\rm\bf\kern0.2em\rule{0.06em}{1.45ex}\kern-0.3em
R$}}

\begin{document}
\title{\large \bf
Normal, cohyponormal and  normaloid weighted composition operators on the Hardy and
weighted Bergman spaces}
 \author{ Mahsa Fatehi and Mahmood Haji Shaabani}

{\maketitle}
\begin{abstract}
If $\psi$ is analytic on the open unit disk $\mathbb{D}$ and $\varphi$ is an analytic self-map
of $\mathbb{D}$, the weighted composition operator
$C_{\psi,\varphi}$ is defined by
$C_{\psi,\varphi}f(z)=\psi(z)f (\varphi (z))$, when $f$ is analytic on $\mathbb{D}$. In this paper, we study normal, cohyponormal, hyponormal and  normaloid weighted composition operators on the Hardy and
weighted Bergman spaces. First, for some weighted Hardy spaces $H^{2}(\beta)$, we prove that if $C_{\psi,\varphi}$ is cohyponormal on $H^{2}(\beta)$, then $\psi$ never vanishes on $\mathbb{D}$ and $\varphi$ is univalent, when $\psi \not \equiv 0$ and $\varphi$ is not a constant function.
Moreover, for $\psi=K_{a}$, where $|a| < 1$, we investigate normal, cohyponormal and hyponormal weighted composition operators $C_{\psi,\varphi}$.
After that, for $\varphi $ which is a hyperbolic or parabolic automorphism, we characterize all normal weighted composition
operators $C_{\psi,\varphi}$, when $\psi \not \equiv 0$ and $\psi$ is analytic on $\overline{\mathbb{D}}$.
 Finally, we find all normal weighted composition operators which are bounded below. \end{abstract}

\footnote{
{\it key words and phrases}:  Weighted Bergman spaces; Hardy space; weighted
composition operator; normaloid operator; cohyponormal operator; normal operator.}

 \section{Introduction}   \par

\bigskip

Let $H(\mathbb{D})$ denote the collection of all holomorphic
functions on the open unit disk $\mathbb{D}$. A function $f$ is called analytic on a closed set $F$ if there exists an open set $U$ such that $f$ is
 analytic on $U$ and $F \subseteq U$. The algebra
$A(\mathbb{D})$ consists of all continuous functions on the
closure of $\mathbb{D}$ that are analytic on ${\mathbb{D}}$.  \par
For $f$ which is analytic on $\mathbb{D}$, we denote by $\hat{f}(n)$ the $n$-th coefficient of $f$ in its Maclaurin series. The Hardy space $H^{2}$ is the collection of all such functions $f$ for which
$$\|f\|_{1}^{2}=\sum_{n=0}^{\infty}|\hat{f}(n)|^{2}<\infty.$$
The space $H^{\infty}$, simply $H^{\infty}$, consists of all functions that are analytic and bounded on ${\mathbb{D}}$.
Recall that for $\alpha >-1$, the weighted Bergman space
$A^{2}_{\alpha}({\mathbb{D}})=A^{2}_{\alpha}$, is the set of
functions $f$ analytic on the unit disk, satisfying the norm
condition
$$\|f\|^{2}_{\alpha}=\int_{\mathbb{D}}|f(z)|^{2}w_{\alpha}(z)dA(z)<\infty,$$
where $w_{\alpha}(z)=(\alpha+1)(1-|z|^{2})^{\alpha}$ and $dA$ is
the normalized area measure. When $\alpha=0$, this gives the
Bergman space $A^{2}(\mathbb{D})=A^{2}$.


Let $e_{w}$ be the linear functional for evaluation at $w$, that is, $e_{w}(f)=f(w)$. Then for functional Hilbert spaces $H$, we let  $K_{w}$ denote the unique function in $H$ which satisfies $\langle f, K_{w} \rangle=f(w)$ for every $f \in H$. In this case, the functional Hilbert space $H$ is called a reproducing kernel Hilbert space.
The weighted Bergman spaces $A_{\alpha}^{2}$ and the Hardy space $H^{2}$ are all reproducing kernel Hilbert
spaces. Let $\gamma=1$ for $H^{2}$ and $\gamma=\alpha+2$ for $A_{\alpha}^{2}$. In $H^{2}$ and $A^{2}_{\alpha}$, we have reproducing kernels $K_{w}(z)=(1-\overline{w}z)^{-\gamma}$ and $\|K_{w}\|_{\gamma}^{2}=(1-|w|^{2})^{-\gamma}$. Moreover,  let $k_{w}$ denote the normalized reproducing kernel given by
$k_{w}(z)=K_{w}(z)/\|K_{w}(z)\|_{\gamma} $.
 \par

Let $\varphi$ be an analytic map of the open unit disk
$\mathbb{D}$ into itself. We define the composition operator $C_{\varphi}$ by $C_{\varphi}(f)=f\circ\varphi$, where $f$ is analytic on $\mathbb{D}$.
If $\psi$ is in $H(\mathbb{D})$ and $\varphi$ is an analytic map
of the unit disk into itself, the weighted composition operator with symbols $\psi$ and $\varphi$ is
the operator $C_{\psi,\varphi}$ which is defined by
$C_{\psi,\varphi}(f)=\psi .(f \circ \varphi )$, where $f$ is analytic on $\mathbb{D}$. If $\psi$ is a bounded analytic
function on $\mathbb{D}$, then the weighted composition operator
$C_{\psi,\varphi}$ is bounded on $H^{2}$ and $A^{2}_{\alpha}$.


 A linear-fractional self-map of $\mathbb{D}$ is a map of the form

\begin{equation}
\varphi(z)=\frac{az+b}{cz+d}
\end{equation}
for some $a,b,c,d \in \mathbb{D}$ such that $ad-bc\neq0$, with the
property that $\varphi(\mathbb{D})\subseteq \mathbb{D}$. We denote the set of
those maps by $\mbox{LFT}(\mathbb{D})$.
It is well-known that the automorphisms of the unit disk, that is,
the one-to-one analytic maps of the unit disk onto itself, are just the
functions
$\varphi(z)=\lambda(a-z)/(1-\bar{a}z)$,
where $|\lambda|=1$ and $|a|<1$ (see, e.g., \cite{c3}). We denote
the class of automorphisms of $\mathbb{D}$ by
$\mbox{Aut}(\mathbb{D})$.

For each $b \in {L^{\infty}(\partial \mathbb{D})}$, we define
the Toeplitz  operator $T_{b}$ on $H^{2}$ by
$T_{b}(f)=P(bf)$, where $P$ denotes the
orthogonal projection of $L^{2}(\partial \mathbb{D})$ onto $H^{2}$.
For each $\psi \in {L^{\infty}(\mathbb{D})}$, we define the
Toeplitz  operator $T_{\psi}$ on $A_{\alpha}^{2}$ by
$T_{\psi}(f)=P_{\alpha}(\psi f)$, where $P_{\alpha}$ denotes the
orthogonal projection of $L^{2}(\mathbb{D},dA_{\alpha})$ onto
$A_{\alpha}^{2}$. Since an orthogonal projection has norm $1$,
clearly $T_{\psi}$ is bounded. If $\psi$ is
a bounded analytic function on $\mathbb{D}$, then the weighted
composition operator can be rewritten as $
C_{\psi,\varphi}=T_{\psi} C_{\varphi}$. \par
 If $\varphi$ is as in Equation (1), then the
adjoint of any linear-fractional composition operator
$C_{\varphi}$, acting on $H^{2}$ and $A_{\alpha}^{2}$, is given by
$C_{\varphi}^{\ast}=T_{g}C_{\sigma}T_{h}^{\ast},$
where
$\sigma(z)=({\overline{a}z-\overline{c}})/({-\overline{b}z+\overline{d}})$
is a self-map of $\mathbb{D}$,
 $g(z)=(-\overline{b}z+\overline{d})^{-\gamma}$,
$h(z)=(cz+d)^{\gamma}$, with $\gamma=1$ for $H^{2}$ and $\gamma=\alpha+2$ for $A^{2}_{\alpha}$. 
From now on,
unless otherwise stated, we assume that $\sigma$, $h$ and $g$ are
given as above. \par

A point $\zeta$ of $\overline{\mathbb{D}}$ is called a fixed point of a self-map $\varphi$ of $\mathbb{D}$ if $\lim_{r \rightarrow 1}\varphi(r \zeta)=\zeta$. We will write $\varphi'(\zeta)$ for $\lim_{r \rightarrow 1}\varphi'(r\zeta)$.
Each analytic self-map $\varphi$ of $\mathbb{D}$ that is neither the identity nor an
elliptic automorphism of $\mathbb{D}$ has associated with it a
unique point $w$ in $\overline{\mathbb{D}}$ that acts like an
attractive fixed point in that $\varphi_{n}(z) \rightarrow w$ as
$n \rightarrow \infty$, where $\varphi_{n}$ denotes $\varphi$
composed with itself $n$ times ($\varphi_{0}$ being the identity
function). The point $w$, called the  Denjoy-Wolff point of $\varphi$, is also characterized as follows:\par $\bullet$ if $|w| < 1$, then $\varphi(w)=w$ and
$|\varphi'(w)|<1$;\par $\bullet$ if $|w| = 1$, then
$\varphi(w)=w$ and $0<\varphi'(w)\leq1$.\\
 More information about Denjoy-Wolff points can be found in \cite[Chapter 2]{cm1} or
\cite[Chapters 4 and 5]{sh}.\par

A map $\varphi \in \mbox{LFT}(\mathbb{D})$ is called parabolic if it has a single fixed point
$\zeta$ in the Riemann sphere $\hat{\mathbb{C}}$ such that $\zeta \in \partial\mathbb{D}$. Let
$\tau(z)=({1+\overline{\zeta}z})/({1-\overline{\zeta}z})$. The map $\tau$ takes the unit disk
onto the right half-plane $\Pi$ and takes $\zeta$ to
$\infty$. The function $\phi=\tau \circ \varphi \circ \tau^{-1}$ is a linear-fractional self-map
of $\Pi$ that fixes only the point $\infty$, so it must have the form $\phi(z)=z+t$ for
some complex number $t$, where $\mbox{Re}(t) \geq 0$.
 Let us call
$t$ the translation number of either $\varphi$ or $\phi$. Note that
if $\mbox{Re}(t)=0$, then $\varphi \in \mbox{Aut}(\mathbb{D})$. Also if
$\mbox{Re}(t)>0$, then $\varphi \not\in \mbox{Aut}(\mathbb{D})$. In \cite[p. 3]{sh}, J. H. Shapiro showed
that among the linear-fractional transformations fixing $\zeta \in
\partial \mathbb{D}$, the parabolic ones are characterized by
$\varphi'(\zeta)=1$.
Let $\varphi
\in \mbox{LFT}(\mathbb{D})$ be parabolic with fixed point $\zeta$ and
translation number $t$. Therefore,
\begin{equation}
\varphi(z)=\frac{(2-t)z+t\zeta}{2+t-t\overline{\zeta}z}.
\end{equation}\par
Recall that an operator $T$ on a Hilbert space $H$ is said to be
normal if $TT^{\ast}=T^{\ast}T$ and  essentially normal if
$TT^{\ast}-T^{\ast}T$ is compact on $H$. Also $T$ is unitary if
$TT^{\ast}=T^{\ast}T=I$. The normal composition operators on
$A^{2}_{\alpha}$ and $H^2$ have symbol $\varphi(z)=az$, where $|a|\leq
1$ (see \cite[Theorem 8.2]{cm1}). An operator $T$ on a Hilbert space $H$ is said to be binormal if $T^{\ast}TTT^{\ast}=TT^{\ast}T^{\ast}T$. It is well known that every normal operator is binormal. In \cite{jkk}, binormal composition operator $C_{\varphi}$ was charactrized when $\varphi$ is a linear-fractional self-map of $\mathbb{D}$.
 If $A^{\ast}A\geq AA^{\ast}$ or, equivalently, $\|Ah\| \geq \|A^{\ast}h\|$ for all vectors $h$, then $A$ is said to be a hyponormal operator.
An operator $A$ is said to be cohyponormal if $A^{\ast}$ is hyponormal. An operator $A$ on a Hilbert space $H$ is normal if and only if
for all vectors $h \in H$, $\|Ah\|=\|A^{\ast}h\|$, so it is not hard to see that an operator $A$ is normal if and only if $A$
and $A^{\ast}$ are hyponormal. Recall that an operator $T$ is said to be normaloid if $\|T\|=r(T)$, where $r(T)$ is the spectral radius of $T$. Then we can see that all normal, hyponormal and cohyponormal operators are normaloid.
 The normal and unitary weighted composition
operators on $H^{2}$ were investigated in \cite{Bourdon} by Bourdon et al.  After that, in \cite{le}, these results were extended to the bigger spaces containing the Hardy and weighted Bergman spaces. Recently,  hyponormal and cohyponormal weighted composition operators have been investigated in \cite{coko} and \cite{mahsa}. In this paper, we work on normal, normaloid, cohyponormal and hyponormal weighted composition operators. In the second section, we extend  \cite[Theorem 3.2]{coko} to some weighted Hardy spaces. In the third section, for $\psi=K_{a}$, where $a \in \mathbb{D}$, we show that if $C_{\psi,\varphi}$ is normal, then $|\varphi(0)|=|a|$.
 In the fourth section,  we state that for $\varphi \in \mbox{Aut}(\mathbb{D})$
and $\psi \not\equiv 0$ which is analytic on
$\overline{\mathbb{D}}$, if $C_{\psi,\varphi}$ is normal on $H^{2}$ or
$A^{2}_{\alpha}$ and $\psi(\zeta)=0$, then $\zeta \in \partial\mathbb{D}$, $\varphi(\zeta)=\zeta$ and
$\zeta$ is not the Denjoy-Wolff point of $\varphi$. Also we prove that if $C_{\psi,\varphi}$  is normal on $H^{2}$ or $A^{2}_{\alpha}$, then $1$ is an eigenvalue of $C_{\psi,\varphi}$, when $\varphi \in \mbox{Aut}(\mathbb{D})$, $\psi$ is analytic on $\overline{\mathbb{D}}$ and $\psi \not \equiv 0$. Furthermore, for $\varphi$ which is a parabolic or hyperbolic automorphism, we give a necessary and sufficient condition for $C_{\psi,\varphi}$ to be normal on $H^{2}$ and $A^{2}_{\alpha}$, when $\psi$ is analytic on  $\overline{\mathbb{D}} $. Finally, we show that for a normal weighted composition operator $C_{\psi,\varphi}$ on a Hilbert space $H$ which contains all the polynomials, $C_{\psi,\varphi}$ is Fredholm if and only if $C_{\psi,\varphi}$ has closed range.\\ \par

\section{Cohyponormal weighted composition operators }
In this section, we provide the generalized result of \cite[Theorem 3.2]{coko} on some weighted Hardy spaces.
We first state the following well-known lemma which was proved in
\cite[p. 1211]{fash} and \cite[p. 1524]{ju}. \\
\par

{\bf Lemma 2.1.} {\it Let $C_{\psi,\varphi}$ be a bounded operator
on $H^{2}$ and $A^{2}_{\alpha}$. For each $w \in \mathbb{D}$,
$C_{\psi,\varphi}^{\ast}K_{w}=\overline{\psi(w)}K_{\varphi(w)}$.}\bigskip

Let $H$ be a Hilbert space. The set of all bounded operators  from $H$ into itself is denoted by $B(H)$.
Now assume that $H$ is a Hilbert space of analytic functions on $\mathbb{D}$. For $f \in H$, let $[f]$ denote the smallest closed subspace of $H$ which contains $\{z^{n}f\}_{n=0}^{\infty}$. If $S \in B(H)$ is the unilateral shift $Sf=zf$, then $[f]$ is the smallest closed subspace of $H$ containing $f$ which is invariant under $S$; moreover, if $[f]=H$, then the function $f$ is called cyclic. Also for $\psi \in H$, we define   a multiplication operator $M_{\psi}:H \rightarrow H$ that for each $f \in H$, $M_{\psi}(f)=\psi f$.  In this section, we assume that $M_{\psi}$ and $S$ are bounded operators, but in general every multiplication operator is not bounded.\\ \par

{\bf Theorem 2.2.} {\it Assume that $H$ is a Hilbert space of analytic functions on $\mathbb{D}$ and the polynomials are dense in $H$. Assume that $\psi  \not \equiv 0$ and $\varphi$ is not a constant function. If $C_{\psi,\varphi}$ is cohyponormal on $H$, then $\psi$ is  cyclic in $H$.}\bigskip \par

{\bf Proof.} Suppose that $C_{\psi,\varphi}$ is cohyponormal. By the Open Mapping Theorem, $\mbox{ker}~C_{\psi,\varphi}=\{0\}$. Then $\mbox{ker}~C_{\psi,\varphi}^{\ast}=\{0\}$. Since $\mbox{ker}~M_{\psi}^{\ast} \subseteq \mbox{ker}~C_{\psi,\varphi}^{\ast}$, we have $\mbox{ker}~M_{\psi}^{\ast}=\{0\}$. \cite[Theorem 2.19, p. 35]{c1} and \cite[Corollary 2.10, p. 10]{c1} imply that $\overline{\mbox{ran}~M_{\psi}}=H$. Then $\psi H$ is dense in $H$. Because the polynomials are dense in $H$, it is easily seen that this is equivalent to saying that polynomial multiples of $\psi$ are dense in $H$, that is, to $\psi$ being a cyclic vector.\hfill $\Box$ \\ \par

Note that by \cite[Corollary 1.5, p. 15]{ha}, if $f \in H^{2}$ is cyclic, then it is an outer function. Then under the conditions of Theorem 2.2, if $C_{\psi,\varphi}$ is cohyponormal on $H^{2}$, then $\varphi$ is an outer function (see \cite[Theorem 3.2]{coko}).\\ \par

{\bf Lemma 2.3.} {\it Let $H$ be a reproducing kernel Hilbert space of analytic functions on $\mathbb{D}$. Assume that for each $w \in \mathbb{D}$,  there is $g \in H$ such that $g(w) \neq 0$. Let $\psi$ be  cyclic in $H$. Then $\psi$ never vanishes on $\mathbb{D}$.}\bigskip \par

{\bf Proof.} Since $\psi$ is cyclic in $H$, $\{p \psi:p ~\mbox{is a polynomial}\}$ is dense in $H$. Let $f \in H$. Then there is a sequence $\{p_{n}\}$ of polynomials  such that $p_{n}\psi \rightarrow f$ as $n \rightarrow \infty$. Suppose that $\psi(w)=0$ for some $w \in \mathbb{D}$. We can see that
$f(w)=\langle f,K_{w}\rangle= \lim_{n \rightarrow \infty }\langle p_{n}\psi,K_{w}\rangle=0$.
Then for each $f \in H$, $f(w)=0$ and it is a contradiction.\hfill $\Box$ \\ \par

Let $H$ be a Hilbert space of analytic functions  on the unit disk. If the monomials $1, z, z^{2},...$ are an orthogonal set of non-zero vectors with dense span in $H$, then $H$ is called a weighted Hardy space. We will assume that the norm satisfies the normalization $\|1\|=1$. The weight sequence for a weighted Hardy space $H$ is defined to be $\beta(n)=\|z^{n}\|$. The weighted Hardy space with weight sequence $\beta(n)$ will be denoted  by $H^{2}(\beta)$.  The inner product on $H^{2}(\beta)$ is given by
$$\langle  \sum_{j=0}^{\infty}a_{j}z^{j},\sum_{j=0}^{\infty}c_{j}z^{j} \rangle=\sum_{j=0}^{\infty}a_{j}\overline{c_{j}}\beta(j)^{2}.$$\par
We require the following corollary, which is a generalization of \cite[Theorem 3.2]{coko}. The proof which shows that $\varphi$ is univalent of the following corollary relies on some ideas from  \cite[Theorem 3.2]{coko}.\\ \par

{\bf Corollary 2.4.} {\it Let $H^{2}(\beta)$ be a weighted Hardy space. Suppose that $\sup \beta(j+1)/\beta(j)$ is finite. Assume that $\psi \not \equiv 0$ and $\varphi$ is not a constant function. If $C_{\psi,\varphi}$ is cohyponormal on $H^{2}(\beta)$, then $\psi$ never vanishes on $\mathbb{D}$ and $\varphi$ is univalent.}\bigskip \par

{\bf Proof.} By Theorem 2.2, Lemma 2.3, \cite[Proposition 2.7]{cm1} and \cite[Theorem 2.10]{cm1}, $\psi$ never vanishes on $\mathbb{D}$.
Assume that there are points $w_{1}$ and $w_{2}$ in $\mathbb{D}$ such that $\varphi(w_{1})=\varphi(w_{2})$ and $w_{1} \neq w_{2}$. Hence
$$C^{\ast}_{\psi,\varphi}(\overline{{\psi}(w_{2})}K_{w_{1}}-\overline{{\psi}(w_{1})}K_{w_{2}})=\overline{{\psi}(w_{2})} \overline{{\psi}(w_{1})}K_{\varphi(w_{1})}-\overline{{\psi}(w_{1})}\overline{{\psi}(w_{2})}K_{\varphi(w_{2})} \equiv 0.$$
We conclude that $0 \in \sigma_{p}(C^{\ast}_{\psi,\varphi})$. Therefore, by \cite[Proposition 4.4, p. 47]{c4}, $0 \in \sigma_{p}(C_{\psi,\varphi})$ and
$C_{\psi,\varphi}(\overline{{\psi}(w_{2})}K_{w_{1}}-\overline{{\psi}(w_{1})}K_{w_{2}})=0.$
Since $\psi$ never vanishes on $\mathbb{D}$, $C_{\varphi}(\overline{{\psi}(w_{2})}K_{w_{1}}-\overline{{\psi}(w_{1})}K_{w_{2}})=0$. Setting $h=K_{w_{2}}/K_{w_{1}}$, we find
$$h\circ \varphi \equiv \overline {\left({\frac{{{\psi}(w_{2})}}{{{\psi}(w_{1})}}}\right)}.$$
Since $\varphi$ is not a constant function, $\varphi(\mathbb{D})$ is an open set by the Open Mapping Theorem. It follows that
$K_{w_{2}}/K_{w_{1}}$ is a constant function and it is a contradiction.\hfill $\Box$ \\ \par

Suppose that $T$ belongs to $B(H^{2})$ or $B(A_{\alpha}^{2})$. Through this paper, the spectrum of $T$, the essential spectrum of $T$ and the point spectrum of $T$ are denoted by $\sigma(T)$, $\sigma_{e}(T)$ and $\sigma_{p}(T)$, respectively.\\ \par

{\bf Remark 2.5.} Suppose that  $C_{\psi,\varphi}$ is cohyponormal on $H^{2}$ or $A^{2}_{\alpha}$ and $\psi \not \equiv 0$.
First, assume that $\varphi$ is not a constant function.
 Since $H^{2}$ and $A_{\alpha}^{2}$ are weighted Hardy spaces,  Theorem 2.2 and Corollary 2.4 imply that $\psi$ is cyclic and $\psi$ never vanishes on $\mathbb{D}$.
Now suppose that $\varphi \equiv c$, where $c$ is a constant number and $|c| < 1$. Assume that there are points $w_{1}$ and $w_{2}$ in $\mathbb{D}$ such that $\psi(w_{1})=0$ and $\psi(w_{2})\neq 0$. From Lemma 2.1, we observe that
 $C^{\ast}_{\psi,\varphi}(K_{w_{1}})=\overline{\psi(w_{1})}K_{c}\equiv 0$.
 Since $C_{\psi,\varphi}$ is cohyponormal and $0 \in \sigma_{p}(C^{\ast}_{\psi,\varphi})$, by \cite[Proposition 4.4, p. 47]{c4}, we have $C_{\psi,\varphi}(K_{w_{1}})=0$. Hence $\psi \cdot K_{w_{1}}\circ\varphi\equiv0$. We conclude that
$\psi(w_{2})K_{w_{1}}(c)=0$
and so $\psi(w_{2})=0$. It is a contradiction. Hence we conclude that $\psi$ never vanishes on $\mathbb{D}$.\\ \par

\section{Normaloid weighted composition operators }

Let $\alpha$ be a complex number of modulus $1$ and $\varphi$ be an analytic self-map of $\mathbb{D}$. Since $\mbox{Re}\left (\frac{\alpha+\varphi(z)}{\alpha-\varphi(z)}\right)$ is a positive harmonic function on $\mathbb{D}$, this function is the poisson integral of a finite positive Borel measure $\mu_{\alpha}$ on $\partial \mathbb{D}$. Let us write $E(\varphi)$ for the closure in $\partial \mathbb{D}$ of the union of the closed supports of the singular parts $\mu_{\alpha}^{s}$ of the measures $\mu_{\alpha}$ as $|\alpha|=1$. In the next Lemma and proposition, the set of points which $\varphi$ makes contact with $\partial \mathbb{D}$ is
$\{\zeta \in \partial \mathbb{D}:\varphi(\zeta)\in \partial \mathbb{D}\}.$ \\ \par

{\bf Lemma 3.1.} \cite[Lemma 3.2]{mahsa} {\it Let $\varphi$ be an analytic self-map of $\mathbb{D}$. Suppose that $\varphi \in A(\mathbb{D})$ and the set of points which $\varphi$ makes contact with $\partial \mathbb{D}$ is finite. Assume that there are a positive integer $n$ and $\zeta \in \partial\mathbb{D}$ such that $E(\varphi_{n})=\{\zeta\}$, where $\zeta$ is the Denjoy-Wolff point of $\varphi$. Let $\psi \in H^{\infty}$ be continuous at $\zeta$. Then
$$r_{\gamma}(C_{\psi,\varphi})=|\psi(\zeta)|\varphi'(\zeta)^{-\gamma/2}.$$
}\bigskip

{\bf Proposition 3.2.} {\it Let $\varphi$ be an analytic self-map of $\mathbb{D}$. Suppose that $\varphi \in A(\mathbb{D})$ and the set of points which $\varphi$ makes contact with $\partial \mathbb{D}$ is finite. Assume that there are a positive integer $n$ and $\zeta \in \partial\mathbb{D}$ such that $E(\varphi_{n})=\{\zeta\}$, where $\zeta$ is the Denjoy-Wolff point of $\varphi$. Suppose that $\psi = K_{a}$ for some $a \in \mathbb{D}$. Let $C_{\psi,\varphi}$ be normaloid on $H^{2}$ or $A_{\alpha}^{2}$. Then
$$\frac{(1-|\varphi(a)|^{2})(1+|a|)}{1-|a|}  \geq \varphi'(\zeta).$$
In particular, if $\varphi'(\zeta)=1$, then $2|a| \geq |\varphi(a)|^{2}(1+|a|)  $.}\bigskip \par

{\bf Proof.} Assume that $C_{\psi,\varphi}$ is normaloid. Let $\gamma =1$ for $H^{2}$ and $\gamma=\alpha+2$ for $A_{\alpha}^{2}$. By  Lemmas 2.1 and 3.1, we can see that
\begin{eqnarray*}
\left|\frac{1}{1-\overline{a}\zeta}\right|^{\gamma}\varphi'(\zeta)^{-\gamma/2}&=&\|C_{\psi,\varphi}\|\nonumber\\
&\geq& \|C_{\psi,\varphi}^{\ast}k_{a}\|_{\gamma}\nonumber\\
&=&|\psi(a)|\frac{\|K_{\varphi(a)}\|_{\gamma}}{\|K_{a}\|_{\gamma}}\nonumber\\
&=&\frac{1}{(1-|a|^{2})^{\gamma}}\left(\frac{1-|a|^{2}}{1-|\varphi(a)|^{2}}\right)^{\gamma/2}.
\end{eqnarray*}
Then
$$\frac{1}{(1-|a|)^{\gamma}}\varphi'(\zeta)^{-\gamma/2} \geq \frac{1}{(1-|a|^{2})^{\gamma}} \left(\frac{1-|a|^{2}}{1-|\varphi(a)|^{2}}\right)^{\gamma/2},$$
 so the result follows. Now suppose that $\varphi'(\zeta)=1$. In this case, after some computation, we can see that $2|a| \geq |\varphi(a)|^{2} +|a| |\varphi(a)|^{2}  $ \hfill $\Box$ \\ \par




{\bf Corollary 3.3.} {\it Let $\varphi$ satisfy the hypotheses of Proposition 3.2. If $C_{\varphi}$ is normaloid on $H^{2}$ or $A_{\alpha}^{2}$, then $1-|\varphi(0)|^{2} \geq \varphi'(\zeta)$. Moreover, if $\varphi'(\zeta)=1$, then $C_{\varphi}$ is not normaloid.}
\bigskip \par

{\bf Proof.} Let $\psi \equiv K_{0}$. By Proposition 3.2, we can see that if $C_{\varphi}$ is normaloid, then $1-|\varphi(0)|^{2} \geq \varphi'(\zeta)$. Now assume that $\varphi'(\zeta)=1$. Suppose that $C_{\varphi}$ is normaloid. Then $1-|\varphi(0)|^{2} \geq 1$. Hence $\varphi(0)=0$ and it is a contradiction.\hfill $\Box$ \\ \par

In Proposition 3.4, we only prove the third part. Proofs of the other parts is similar to part (c) and follows from the definitions of hyponormal and cohyponormal operators.\\ \par

{\bf Proposition 3.4.} {\it Let $\varphi$ be an analytic self-map of $\mathbb{D}$ and $\psi=K_{a}$ for some $a \in \mathbb{D}$. The following statements hold on  $H^{2}$ or $A_{\alpha}^{2}$.\\
 (a) If $C_{\psi,\varphi}$ is cohyponormal, then $|\varphi(0)| \geq |a|$. \\
 (b) If $C_{\psi,\varphi}$ is hyponormal, then $|\varphi(0)| \leq |a|$. \\
 (c) If $C_{\psi,\varphi}$ is normal, then $|\varphi(0)| = |a|$.}\bigskip \par

{\bf Proof.} (c) Since $K_{0} \equiv 1$, by \cite[Proposition 2.16, p. 34]{c1} and Lemma 2.1, wee see that
\begin{equation}
\|\overline{\psi(0)}K_{\varphi(0)}\|_{\gamma}=\|C_{\psi,\varphi}^{\ast}K_{0}\|_{\gamma}=\|C_{\psi,\varphi}K_{0}\|_{\gamma}=\|\psi\|_{\gamma}=\|K_{a}\|_{\gamma},
\end{equation}
where $\gamma=1$ for $H^{2}$ and $\gamma=\alpha+2$ for $A_{\alpha}^{2}$.
We know that $\psi=K_{a}$, so Equation (3)
 shows that $|\varphi(0)|=|a|$.\hfill $\Box$ \\ \par


Let $\varphi \in
\mbox{LFT}(\mathbb{D})$. It is easy to see that $\varphi$ must belong to one of the following three disjoint classes:\\
$\bullet$ Automorphism of $\mathbb{D}$.\\
$\bullet$ Non-automorphism of $\mathbb{D}$ with $\overline{\varphi(\mathbb{D})} \subseteq \mathbb{D}$.\\
 $\bullet$  Non-automorphism of $\mathbb{D}$ with $\varphi(\zeta)=\eta$ for some $\zeta,\eta \in \partial \mathbb{D}$. \par
 Let $\varphi \in \mbox{LFT}(\mathbb{D})$ such that $\overline{\varphi(\mathbb{D})
} \subseteq \mathbb{D}$. Then by \cite[Theorem 2.48]{cm1}, $\varphi$ has a fixed point $p \in \mathbb{D}$.
Suppose that $\varphi \in \mbox{LFT}(\mathbb{D})$ such that $\overline{\varphi(\mathbb{D})
} \subseteq \mathbb{D}$ or $\varphi$  is the identity or an automorphism of $\mathbb{D}$ with a fixed point in $\mathbb{D}$. All normal weighted composition operators $C_{\psi,\varphi}$ of these types were found (see  \cite[Theorem 10]{Bourdon}, \cite[Theorem 3.7]{coko} and \cite[Theorem 4.3]{le}).
Also suppose $\varphi \in \mbox{Aut}(\mathbb{D})$ which has no fixed point in  $\mathbb{D}$ and $\psi= K_{a}$ for some $a \in \mathbb{D}$; all normal weighted composition operators $C_{\psi,\varphi}$ on $H^{2}$ and $A^{2}_{\alpha}$ are characterized in Theorem 4.5. Bourdon et al. in \cite[Proposition 12]{Bourdon} obtained a condition that reveals what is required for normality of a weighted composition operator $C_{\psi,\varphi}$ on $H^{2}$, where $\varphi$ is a linear-fractional and $\psi=K_{\sigma(0)}$ (also by the similar proof, an analogue result holds on $A_{\alpha}^{2}$). In the following corollary, for $\psi=K_{a}$ and $\varphi \in
\mbox{LFT}(\mathbb{D})$,  where $\varphi(\zeta)=\eta$ for some $\zeta,\eta \in \partial \mathbb{D}$ and $a \in \mathbb{D}$, we investigate normal weighted composition operators $C_{\psi,\varphi}$ on $H^{2}$ and $A_{\alpha}^{2}$.\\ \par

{\bf Corollary 3.5.} {\it Suppose that $\varphi \in
\mbox{LFT}(\mathbb{D})$ is not an automorphism and $\varphi(\zeta)=\eta$ for some $\zeta,\eta \in \partial \mathbb{D}$. Assume that $\psi=K_{a}$ for some $a \in \mathbb{D}$. If $C_{\psi,\varphi}$ is normal on $H^{2}$ or $A_{\alpha}^{2}$, then $\varphi$ is a parabolic non-automorphism and $|a|=|t/(2+t)|=|\sigma(0)|$, where $t$ is the translation number of $\varphi$.}\bigskip \par

{\bf Proof.}  Let $C_{\psi,\varphi}$ be normal on $H^{2}$ or $A^{2}_{\alpha}$. Then $C_{\psi,\varphi}$ is essentially normal. Since $\psi$ never vanishes on $\partial\mathbb{D}$, we conclude from \cite[Theorem 2.6]{fash} and  \cite[Remark 2.7]{fash} that $\varphi$ is  a parabolic non-automorphism and the result follows from Proposition 3.4 and Equation (2).\hfill $\Box$ \\ \par

\section{Normal weighted composition operators}


Each disk automorphism $\varphi$ such that $\varphi \neq I$ has two fixed points on the sphere, counting multiplicity. The automorphisms are classified according to the location of their fixed points: elliptic if one  fixed point is in $\mathbb{D}$ and a second fixed point is in the complement of the closed disk, hyperbolic if both fixed points are in $\partial \mathbb{D}$, and parabolic if there is one fixed point in $\partial \mathbb{D}$ of multiplicity two (see \cite{cm1} and \cite{sh}). Let $\varphi$ be an automorphism of $\mathbb{D}$. In  \cite{fash} and \cite{fk1}, the present authors  investigated essentially normal weighted composition operator $C_{\psi,\varphi}$, when $\psi \in A(\mathbb{D})$ and $\psi(z) \neq 0$ for each $z \in \overline{\mathbb{D}}$.  In this section, we just assume that $\psi$ is analytic on $\overline{\mathbb{D}}$ and we attempt to find all normal weighted composition operators $C_{\psi,\varphi}$. Also we will show that $\psi$ never vanishes on $\overline{\mathbb{D}}$.\\ \par

{\bf Lemma 4.1.} {\it Let $\varphi \in \mbox{Aut}(\mathbb{D})$ and $f \in A(\mathbb{D})$. Then $T_{f}^{\ast}C_{\varphi}-C_{\varphi}T_{f\circ \varphi^{-1}}^{\ast}$ is compact on $H^{2}$ and $A^{2}_{\alpha}$.}\bigskip \par

{\bf Proof.} We know that $\sigma=\varphi^{-1}$. It is not hard to see that
$$C_{\varphi}^{\ast}T_{f}=T_{g}C_{\sigma}T_{h}^{\ast}T_{f}=C_{\sigma}T_{g\circ \varphi}T_{h}^{\ast}T_{f}$$
and

$$T_{f\circ \varphi^{-1}}C_{\varphi}^{\ast}=T_{f\circ \varphi^{-1}}C_{\sigma}T_{g \circ \varphi}T_{h}^{\ast}=C_{\sigma}T_{f}T_{g \circ \varphi}T_{h}^{\ast}.$$
 Since $C_{\varphi}C_{\sigma}=I$, by \cite[Proposition 7.22]{do} and \cite[Corollary 1(c)]{ot}, $C_{\varphi}^{\ast}T_{f}-T_{f \circ \varphi^{-1}}C_{\varphi}^{\ast}$ is compact and the result follows.  \hfill $\Box$ \\ \par

In this section, we assume that $\varphi(z)=\lambda(a-z)/(1-\overline{a}z)$ and $w(z)=(1-\overline{a}z)^{\gamma}\psi(z)$, where $a \in \mathbb{D}$, $|\lambda|=1$, $\psi \in A(\mathbb{D})$ and $\gamma=1$ for $H^{2}$ and $\gamma=\alpha+2$ for $A^{2}_{\alpha}$.
 We will use the notation $A \equiv B$
to indicate that the difference of
the two bounded operators $A$ and $B$ is compact.\\ \par

{\bf Proposition 4.2.} {\it Let $\psi \in A(\mathbb{D})$ and  $\varphi \in \mbox{Aut}(\mathbb{D})$. If $C_{\psi,\varphi}$ is hyponormal on $H^{2}$ or  $A^{2}_{\alpha}$, then for each $\zeta \in \partial \mathbb{D}$, $|w(\zeta)|-|w(\varphi(\zeta))| \geq 0$. Moreover, if $C_{\psi,\varphi}$ is cohyponormal on $H^{2}$ or  $A^{2}_{\alpha}$, then for each $\zeta \in \partial \mathbb{D}$, $|w(\varphi(\zeta))| -|w(\zeta)| \geq 0$.}\bigskip \par

{\bf Proof.} Suppose that $C_{\psi,\varphi}$ is hyponormal. By  \cite[Proposition 7.22]{do}, \cite[Corollary 1(c)]{ot} and the preceding lemma, we can see that
\begin{eqnarray*}
C_{\psi,\varphi}C_{\psi,\varphi}^{\ast}&=&T_{\psi}C_{\varphi}T_{g}C_{\sigma}T_{h}^{\ast}T_{\psi}^{\ast}\nonumber\\
&=&T_{\psi}C_{\varphi}C_{\sigma}T_{g \circ \varphi}T_{h}^{\ast}T_{\psi}^{\ast}\nonumber\\
& \equiv &T_{|\psi|^{2} \cdot \overline{h} \cdot g \circ \varphi}
\end{eqnarray*}
and
\begin{eqnarray*}
C_{\psi,\varphi}^{\ast}C_{\psi,\varphi}&=&T_{g}C_{\sigma}T_{h}^{\ast}T_{\psi}^{\ast}T_{\psi}C_{\varphi}\nonumber\\
&=&T_{g}C_{\sigma}T_{h\psi}^{\ast}C_{\varphi}T_{\psi \circ \varphi^{-1}}\nonumber\\
& \equiv &T_{g}C_{\sigma}C_{\varphi}T_{(h\psi)\circ\varphi^{-1}}^{\ast}T_{\psi \circ \varphi^{-1}}\nonumber\\
& \equiv &T_{|\psi \circ \varphi^{-1}|^{2} \cdot g \cdot \overline{h} \circ \varphi^{-1}}.
\end{eqnarray*}
Let $\varphi(a)=0$ for $a \in \mathbb{D}$. After some computation, we see that
$\overline{h(z)}g(\varphi(z))=|1-\overline{a}z|^{2\gamma}/(1-|a|^{2})^{\gamma},$
where $\gamma=1$ for $H^{2}$ and $\gamma=\alpha+2$ for $A_{\alpha}^{2}$. Hence by \cite[Corollary 2.6]{e} and \cite[Corollary 1.3]{cob}, we have
\begin{eqnarray*}
\sigma_{e}(C_{\psi,\varphi}^{\ast}C_{\psi,\varphi}-C_{\psi,\varphi}C_{\psi,\varphi}^{\ast})&=& \frac{1}{(1-|a|^{2})^{\gamma}}\sigma_{e}(T_{|w \circ \varphi^{-1}|^{2}-|w|^{2}})\nonumber\\
&=&\left\{\frac{|w(\varphi^{-1}(\zeta))|^{2}-|w(\zeta)|^{2}}{(1-|a|^{2})^{\gamma}}:\zeta \in \partial \mathbb{D}\right\}.
\end{eqnarray*}
Since $C_{\psi,\varphi}^{\ast}C_{\psi,\varphi} \geq C_{\psi,\varphi}C_{\psi,\varphi}^{\ast}$, $|w(\varphi^{-1}(\zeta))|^{2}-|w(\zeta)|^{2} \geq 0$ for each $\zeta \in \partial \mathbb{D}$ and so $|w(\zeta)|^{2}-|w(\varphi(\zeta))|^{2} \geq 0$ for any $\zeta \in \partial \mathbb{D}$. Therefore, the conclusion follows. The idea of the proof of the result for cohyponormal operator $C_{\psi,\varphi}$ is similar to hyponormal operator, so it is left for the reader. \hfill $\Box$ \\ \par

Suppose $\psi$ is a non-constant analytic function on $\overline{\mathbb{D}}$ and $\psi$ never vanishes on $\mathbb{D}$.  By \cite[Exercise 1, p. 129]{c3}, we see that $|\psi|$ assumes its minimum value on $\partial \mathbb{D}$. Assume that $C_{\psi,\varphi}$ is essentially normal on $H^{2}$ or $A^{2}_{\alpha}$. Now let $\psi(\zeta)=0$ for some $\zeta \in \partial \mathbb{D}$. In the following proposition, we show that $\varphi(\zeta)=\zeta$ and $\zeta$ is not the Denjoy-Wolff point of $\varphi$. \\ \par

{\bf Proposition 4.3.} {\it Let $\varphi \in
\mbox{Aut}(\mathbb{D})$. Suppose that $\psi$ is analytic on
$\overline{\mathbb{D}}$ and $\psi$ never vanishes on $\mathbb{D}$. Let
$C_{\psi,\varphi}$ be essentially normal on $H^{2}$ or $A^{2}_{\alpha}$. If
$\psi(\zeta)=0$ for some $\zeta \in  \partial \mathbb{D}$, then
 $\varphi(\zeta)=\zeta$ and
$\zeta$ is not the Denjoy-Wolff point of $\varphi$.}\bigskip \par

{\bf Proof.} Let $C_{\psi,\varphi}$ be essentially normal on $H^{2}$ or
$A^{2}_{\alpha}$.
Assume that there is $\zeta \in
\partial\mathbb{D}$ such that $\psi(\zeta)=0$. By  \cite[Theorem 3.2]{fash} and \cite[Theorem 3.3]{fk1}, $\psi$ is zero on $B=\{\zeta,
\varphi(\zeta), \varphi_{2}(\zeta),...\}$. It is easy to see
that $B \subseteq \partial\mathbb{D}$. If $|B|=\infty$, then
$\psi\equiv0$ and it is a contradiction. Let $|B|=N$ with
$N<\infty$. We have $\varphi_{N} \in \mbox{Aut}(\mathbb{D})$ and it is not hard to see that
for each $b \in B$, ${\varphi}_{N}(b)=b$. Since every automorphism
of $\mathbb{D}$ has at most two fixed points, we conclude that $N\leq2$. Hence $B=\{\zeta\}$ or
$B=\{\zeta, \varphi(\zeta)\}$. We claim that $\zeta$ is the fixed point of $\varphi$. This may be seen as follows. If $B$ has only one element, then $\varphi(\zeta)=\zeta$. Now let $N=2$. One can easily see that if $c$ is a
fixed point of $\varphi$, then $c$ is also the fixed point of
$\varphi_{2}$. Therefore, the set of all fixed points of $\varphi$ is a
subset of $B$. In this case, assume that
$\varphi(\zeta)$ is the fixed point of $\varphi$ and so $\varphi(\varphi(\zeta))=\varphi(\zeta)$. Since $\zeta \in B$ and $|B|=2$, we conclude that
$\varphi_{2}(\zeta)=\zeta$. Therefore, $\varphi(\zeta)=\zeta$. Hence $B=\{\zeta\}$ and $\zeta$ is the fixed point of $\varphi$.
 Since $|w|=|w \circ \varphi|$ on $\partial \mathbb{D}$, by  \cite[Exercise 6, p. 130]{c3}, there is a constant $\lambda$, $|\lambda|=1$ such that $w \circ \varphi=\lambda w$.
 Suppose that $\zeta$
is the Denjoy-Wolff point of $\varphi$. Since for each $z \in
\mathbb{D}$, $\{z\}$ is a compact set, the Denjoy-Wolff Theorem
implies that
$$0=|w(\zeta)|=\lim_{n\rightarrow\infty} |w(\varphi_{n}(z))|=\lim_{n\rightarrow\infty}|\lambda^{n}w(z)|=|w(z)|.$$
Hence $\psi\equiv0$ and it is a contradiction.
 \hfill $\Box$ \\ \par

If $C_{\psi,\varphi}$ is normal on $H^{2}$ or $A^{2}_{\alpha}$, then by Remark 2.5, $\psi$ never vanishes on $\mathbb{D}$.
In Proposition 4.3, we saw that for $\varphi \in \mbox{Aut}(\mathbb{D})$ which is not a hyperbolic automorphism and $\psi$ that is analytic on $\overline{\mathbb{D}}$, if $C_{\psi,\varphi}$ is normal on $H^{2}$ or $A^{2}_{\alpha}$, then $\psi$ never vanishes on $\overline{\mathbb{D}}$. \\ \par

{\bf Proposition 4.4.} {\it Let $\varphi \in \mbox{Aut}(\mathbb{D})$ and $\varphi(a)=0$  for some $a \in \mathbb{D}$. Assume that
$\psi$ is analytic on $\overline{\mathbb{D}}$ and $\psi \not  \equiv 0$. If  $C_{\psi,\varphi}$ is  normal on $H^{2}$ or  $A^{2}_{\alpha}$, then $w$ is an eigenvector
for the operator $C_{\varphi}$  and the corresponding
$C_{\varphi}$-eigenvalue for $w$ is 1. }\bigskip \par

{\bf Proof.} Suppose that $C_{\psi,\varphi}$ is normal. Proposition 4.2,  implies that $|w|=|w\circ\varphi|$ on $\partial\mathbb{D}$. Since $\psi$ never vanishes on $\mathbb{D}$, by \cite[Exercise 6, p. 130]{c3}, we conclude that $C_{\varphi}(w)=\lambda w$, where $|\lambda|=1$. If $\varphi$ is an elliptic automorphism with a fixed point $t \in \mathbb{D}$, then $w(t)=w(\varphi(t))=\lambda w(t)$. Since $\psi$ never vanishes on $\mathbb{D}$, $\lambda=1$. Now suppose that $\varphi$ is a parabolic or hyperbolic automorphism with a Denjoy-Wolff point $\zeta$. Then
$$w(\zeta)=\lim _{r \rightarrow 1}w(\varphi(r \zeta))=\lambda \lim_{r \rightarrow 1}w(r \zeta)=\lambda w(\zeta).$$
By Proposition 4.3, $\psi(\zeta) \neq 0$ and so $\lambda=1$. Therefore, we conclude that $C_{\varphi}(w)=w$. \hfill $\Box$ \\ \par

Let $\varphi$ be an elliptic automorphism or the identity. As we stated before Corollary 3.5, all normal weighted composition operators $C_{\psi,\varphi}$ of these types were found. Also we must say that in this case $\psi$ never vanishes on $\overline{\mathbb{D}}$. In the next theorem, for $\varphi$, not the identity and not an elliptic automorphism of $\mathbb{D}$, which is in $\mbox{Aut}(\mathbb{D})$, we show that constant multiples of $K_{\sigma(0)}$ are the only examples for $\psi$ that $C_{\psi,\varphi}$ are normal, where $\psi$ is analytic on
$\overline{\mathbb{D}}$. It is interesting that again $\psi$ never vanishes on $\overline{\mathbb{D}}$ and these weighted composition operators are actually a constant multiples of unitary weighted composition operators (see \cite[Theorem 6]{Bourdon} and \cite[Corollary 3.6]{le}).\\ \par

 {\bf Theorem 4.5.} {\it Assume that $\varphi$, not the identity and not an elliptic automorphism of $\mathbb{D}$, is in $\mbox{Aut}(\mathbb{D})$. Suppose that $\psi$ is analytic on
$\overline{\mathbb{D}}$ and $\psi \not\equiv 0$. Then $C_{\psi,\varphi}$ is normal on  $H^{2}$ or $A^{2}_{\alpha}$ if and only if $\psi=\psi(0)K_{\sigma(0)}$; hence $\psi$ never vanishes on $\overline{\mathbb{D}}$.}\bigskip \par

{\bf Proof.} Suppose that $C_{\psi,\varphi}$ is normal and $\varphi(a)=0$. By Proposition 4.4, $w \circ \varphi=w$ on $ \mathbb{D}$. Then for each integer $n$, $\psi(z)(1-\overline{a}z)^{\gamma}=\psi(\varphi_{n}(z))(1-\overline{a}\varphi_{n}(z))^{\gamma}$ on $\mathbb{D}$, where $\gamma= 1$ for $H^{2}$ and $\gamma =\alpha+2$ for $A_{\alpha}^{2}$. Let $\zeta$ be a Denjoy-Wolff point of $\varphi$. Since for each $z \in \mathbb{D}$, $\{z\}$ is a compact set, we have
$$\psi(\zeta)(1-\overline{a}\zeta)^{\gamma}= \lim_{n \rightarrow \infty}\psi(\varphi_{n}(z))(1-\overline{a}\varphi_{n}(z))^{\gamma}=\psi(z)(1-\overline{a}z)^{\gamma}.$$
Hence $\psi(z)=\frac{ \psi(\zeta)(1-\overline{a}\zeta)^{\gamma}}{(1-\overline{a}z)^{\gamma}}$ and so
 $\psi=\psi(0)K_{\sigma(0)}$.\\
Conversely, it is not hard to see that the fact that for a constant number $c$, $C_{c\psi,\varphi}$ is normal implies that $C_{\psi,\varphi}$ is also. Then without loss of generality, we assume that $\psi= K_{\sigma(0)}$, where $\gamma=1$ for $H^{2}$ and $\gamma =\alpha+2$ for $A^{2}_{\alpha}$. Observe that $g \circ \varphi=\|K_{\sigma(0)}\|^{2}_{\gamma}/K_{\sigma(0)}$ and $h \psi \equiv 1$. Since $\sigma= \varphi^{-1}$, we get $$C_{\psi,\varphi}^{\ast}C_{\psi,\varphi}=T_{g}C_{\sigma}T_{h \psi}^{\ast}T_{\psi}C_{\varphi}=T_{g \cdot \psi \circ \varphi^{-1}}$$
and
$$C_{\psi,\varphi}C_{\psi,\varphi}^{\ast}=T_{\psi}C_{\varphi}T_{g}C_{\sigma}T_{h \psi}^{\ast}=T_{\psi \cdot g \circ \varphi}.$$
After some computation, one can see that $\psi \cdot g \circ \varphi=g \cdot \psi \circ \varphi^{-1}=\|K_{\sigma(0)}\|_{\gamma}^{2}$.
Hence $C_{\psi,\varphi}$ is normal.
 \hfill $\Box$ \\ \par

{\bf Lemma 4.6.} {\it Assume that $\varphi \in \mbox{Aut}(\mathbb{D})$. Suppose that  for some $\zeta \in \partial \mathbb{D}$, $\{ \varphi_{n}(\zeta):n~ \mbox{is a positive integer}\}$ is a finite set. If $\varphi$ is parabolic or hyperbolic, then $\zeta$ is the fixed point of $\varphi$.}\bigskip \par

{\bf Proof.} Let $\{\zeta,\varphi(\zeta),\varphi_{2}(\zeta),...\}$ be a finite set. Then there is an integer $N$ such that $\varphi_{N}(\zeta)=\zeta$, so $\zeta$ is the fixed point of $\varphi_{N}$. It is not hard to see that $\varphi$ is parabolic or hyperbolic if and only if $\varphi_{N}$ is parabolic or hyperbolic, respectively and fixed points of $\varphi$ and $\varphi_{N}$ are the same. Hence $\zeta$ is the fixed point of $\varphi$.\hfill $\Box$ \\ \par

Let $\psi \in A(\mathbb{D})$. If $C_{\psi,\varphi}$ is cohyponormal on $H^{2}$ or $A^{2}_{\alpha}$, then $\psi$ never vanishes on $\mathbb{D}$ or $\psi \equiv 0$
 (see Remark 2.5). Therefore, by Maximum Modulus Theorem and \cite[Exercise 1, p. 129]{c3}, there are $\zeta_{1},\zeta_{2} \in \partial \mathbb{D}$, with $|w(\zeta_{1})| \leq |w(z)|$ and $|w(z)| \leq |w(\zeta_{2})|$  for all $z \in \mathbb{D}$. In the following theorem, we assume that $\zeta_{1}$ and $\zeta_{2}$ are given as above. \\ \par

 {\bf Theorem 4.7.} {\it Suppose that  $\psi \in A(\mathbb{D})$. Let $\{\eta \in \partial \mathbb{D}: |w(\eta)|=|w(\zeta_{1})|\}$ and $\{\eta \in \partial \mathbb{D}: |w(\eta)|=|w(\zeta_{2})|\}$ be finite sets. Suppose that $C_{\psi,\varphi}$ is cohyponormal on $H^{2}$ or $A^{2}_{\alpha}$.  The following statements hold.\\
 (a) If $\varphi$ is a parabolic automorphism, then $|w|$ is a constant function on $\partial \mathbb{D}$. Moreover, if $\psi \not \equiv 0$, then $C_{\psi,\varphi}$ is normal and $\psi=\psi(0)K_{\sigma(0)}$.\\
 (b) If $\varphi$ is a hyperbolic automorphism, then $\zeta_{1}$ and $\zeta_{2}$ are the fixed points of $\varphi$.\\
 }\bigskip \par

  {\bf Proof.} (a) Suppose that $\varphi$ is a parabolic automorphism. Since $C_{\psi,\varphi}$ is cohyponormal, by Proposition 4.2, for each positive integer $n$, we have
 $$|w(\varphi_{n}(\zeta_{2}))| \geq |w(\varphi_{n-1}(\zeta_{2}))| \geq |w(\varphi_{n-2}(\zeta_{2}))| \geq ...\geq |w(\zeta_{2})|$$
 and
 $$|w(\zeta_{1})| \geq |w(\varphi^{-1}(\zeta_{1}))| \geq |w(\varphi^{-1}_{2}(\zeta_{1}))| \geq ...\geq |w(\varphi^{-1}_{n}(\zeta_{1}))|.$$
 It is not hard to see that $\varphi^{-1}$ is parabolic and $\varphi$ and $\varphi^{-1}$ have the same fixed point.
 Then by Lemma 4.6 and the statement which was stated before Theorem 4.7, we can see that $\zeta_{1}=\zeta_{2}$. Then $|w|$ is constant on  $\partial \mathbb{D}$. Now assume that $\psi \not \equiv 0$. We know that $\psi$ never vanishes on $\mathbb{D}$ (see Corollary 2.4), so by \cite[Exercise 2, p. 129]{c3}, $w$ is a constant function. Then $\psi=\psi(0)K_{a}$ and by Theorem 4.5, $C_{\psi,\varphi}$ is normal.\\
 (b)
Suppose that $\varphi$ is a hyperbolic automorphism. By the proof of part (a) and Lemma 4.6, we can see that $\zeta_{1}$ and $\zeta_{2}$
are the fixed points of $\varphi$.\\
 \hfill $\Box$ \\ \par



For $\psi \in H^{\infty}$ and $\varphi$ which is an elliptic automorphism of $\mathbb{D}$,  cohyponormality and normality of a weighted composition operator $C_{\psi,\varphi}$ on $H^{2}$ are equivalent (see \cite[Proposition 3.17]{coko}). In the previous theorem, we showed that if $\varphi$ is a parabolic automorphism and $\psi$ and $w$ satisfy the hypotheses of this theorem, then $C_{\psi,\varphi}$ is cohyponormal on $H^{2}$ or $A^{2}_{\alpha}$ if and only if $C_{\psi,\varphi}$ is normal and we saw that $\psi=\psi(0)K_{\sigma(0)}$.\par
Recall that a bounded linear  operator $T$ between two Banach spaces is Fredholm if it is invertible modulo compact operators.  We say that an operator $A \in B(H)$ is bounded below if there is a constant $c > 0$ such that
$c\| h \| \leq  \|A(h)\|$
for all $h \in H$.
Moreover, we know that a normal operator $N$ on a Hilbert space $H$ is bounded below if and only if $N$ is invertible (see \cite[Exercise 15, p. 36]{c1}). By this fact, the statements (a) and (c) in Theorem 4.8 are equivalent.   \\ \par

{\bf Theorem 4.8.} {\it Suppose that $C_{\psi,\varphi}$ is a normal operator on a Hilbert space $H$ of analytic functions on $\mathbb{D}$. Assume that all the polynomials belong to $H$. The following statements are equivalent.\\
(a) The operator $C_{\psi,\varphi}$ is bounded below.\\
(b) The operator $C_{\psi,\varphi}$ is Fredholm.\\
(c) The operator $C_{\psi,\varphi}$ is invertible.}\bigskip \par

{\bf Proof.}
(b) implies (c). Suppose that $C_{\psi,\varphi}$ is Fredholm. Since by \cite[Corollary 2.4, p. 352]{c1}, $\mbox{dim}~\mbox{ker}~C_{\psi,\varphi} < \infty$, it is not hard to see that $\psi \not \equiv 0$. We claim that $\varphi$ is not a constant function. Assume that $\varphi \equiv c$, where $|c| < 1$. It is not hard to see that for each $n$, $z^{n}(z-c) \in \mbox{ker}~C_{\psi,\varphi}$, so  $\mbox{dim}~\mbox{ker}~C_{\psi,\varphi}=\infty $ and it is a contradiction. By the Open Mapping Theorem, we can see that $0 \not\in \sigma_{p}(C_{\psi,\varphi})$. Assume that $C_{\psi,\varphi}$ is not invertible. Then by \cite[Proposition 4.6, p. 359]{c1}, $0 \in \sigma_{p}(C_{\psi,\varphi})$ and it is a contradiction.\par
(c) implies (b). This is clear.\hfill $\Box$ \\

Assume $\psi \not\equiv 0$ and $\varphi$ is not a constant function. By the Open Mapping Theorem, it is clear that $\mbox{ker}~C_{\psi,\varphi}=(0)$. Then $C_{\psi,\varphi}$ has closed range if and only if $C_{\psi,\varphi}$ is bounded below. We know that $C_{\psi,\varphi}$ on $H^{2}$ or $A_{\alpha}^{2}$ is invertible if and only if $\varphi \in \mbox{Aut}(\mathbb{D})$ and $\psi \in H^{\infty}$ is bounded away from zero on $\mathbb{D}$ (see \cite[Theorem 3.4]{Bourdon2}). As we stated before if $\varphi$ is an elliptic automorphism or the identity, all normal weighted composition operators were found; moreover others were characterized in Theorem 4.8. Then closed range weighted composition operators  on $H^{2}$ and  $A_{\alpha}^{2}$ which are normal were found.

\bigskip
{
M. Fatehi, Department of Mathematics, Shiraz Branch, Islamic Azad
University, Shiraz, Iran. \par E-mail: fatehimahsa@yahoo.com \par
M. Haji Shaabani, Department of Mathematics, Shiraz University of Technology, P. O. Box 71555-313,
Shiraz, Iran.\par E-mail: shaabani@sutech.ac.ir\par}

\end{document}